\documentclass{amsart}
\usepackage{amsmath,amssymb}
\usepackage[curve,matrix,arrow]{xy}
\usepackage{pb-xy}
\usepackage{pb-diagram}
\allowdisplaybreaks
\newtheorem{thm}{Theorem}
\newtheorem{lem}[thm]{Lemma}
\newtheorem{defn}[thm]{Definition}
\newtheorem{prop}[thm]{Proposition}

\newtheorem{cor}[thm]{Corollary}

\newtheorem{expls}[thm]{Examples}
\newtheorem{rem}[thm]{Remark}

\makeatletter
\def\proofname{{\it Proof}\,:}
\def\aproofname{{\it Proof of }}
\def\pproofname{{\it. }}
\def\@proof[#1]{\aproofname{\it #1}\pproofname}
\renewenvironment{proof}{\par\noindent\@ifnextchar[{\@proof}{{\proofname}\quad }}{{\unskip\nobreak\hfill{\it\qedsymbol}}\par\vskip 9pt}

\ifx\undefined\bysame\newcommand{\bysame}{\leavevmode\hbox to3em{\hrulefill}\,}\fi
\makeatother
\numberwithin{equation}{section}
\numberwithin{thm}{section}
\newcommand{\Loop}[1]{{\Omega}{(#1)}}

\newcommand{\product}{\mathrm{\Pi}}

\newcommand{\fatvee}{\operatorname{T}}

\newcommand{\ad}{\operatorname{ad}}
\newcommand{\id}{\operatorname{id}}
\newcommand{\incl}{\mbox{in}}
\newcommand{\proj}{\mbox{pr}}
\newcommand{\comp}{\circ}

\newcommand{\midvert}{\,\mathstrut\vrule\,}
\newcommand{\integral}{\mathbb Z}

\newcommand{\complex}{\mathbb C}

\newcommand{\gtwo}[1]{\mathrm{G}_{2}}
\newcommand{\ffour}[1]{\mathrm{F}_{4}}
\newcommand{\esix}[1]{\mathrm{E}_{6}}
\newcommand{\eseven}[1]{\mathrm{E}_{7}}
\newcommand{\eeight}[1]{\mathrm{E}_{8}}
\newcommand{\Min}{\operatorname{Min}}
\newcommand{\Max}{\operatorname{Max}}
\gdef\empty{}
\def\TC#1{\operatorname{\ifx\args\empty\operatorname{\mathcal{TC}}\else\mathcal{TC}}(#1)\fi}

\def\TCM#1{\def\args{#1}\ifx\args\empty\operatorname{\mathcal{TC^{M}}}\else\operatorname{\mathcal{TC}^{M}}(#1)\fi}
\def\TCcuplen#1#2{\operatorname{\mathcal{Z}_{#2}}(#1)}
\def\TCwgt#1#2{\operatorname{wgt_{#2}}(#1)}

\def\product{\mathrm{\Pi}}

\def\fatvee{\operatorname{T}}

\def\double#1{\operatorname{\it d}(#1)}
\def\secat#1{\operatorname{secat}(#1)}
\def\genus#1{\operatorname{Genus}(#1)}
\def\Category#1{\underline{\underline{\mathcal{#1}}}}
\newcommand{\Path}[1]{\operatorname{\mathcal P}(#1)}

\newcommand{\PathB}[1]{\operatorname{{\mathcal P}_{B}}(#1)}
\newcommand{\fLoop}[1]{\operatorname{{\mathcal L}}(#1)}
\newcommand{\Track}[1]{\def\args{#1}\ifx\args\empty\operatorname{\mathcal{I}}\else\operatorname{\mathcal{I}}(#1)\fi}
\newcommand{\TrackB}[1]{\def\args{#1}\ifx\args\empty\operatorname{\mathcal{I}}_{B}\else\operatorname{\mathcal{I}}_{B}(#1)\fi}
\newcommand{\TrackBB}[1]{\def\args{#1}\ifx\args\empty\operatorname{\mathcal{I}}^{B}_{B}\else\operatorname{\mathcal{I}}^{B}_{B}(#1)\fi}

\newcommand{\simBB}{\simeq^{B}_{B}}
\newcommand{\Susp}[1]{\def\args{#1}\ifx\args\empty\operatorname{\mathcal{S}}\else\operatorname{\mathcal{S}}(#1)\fi}
\renewcommand{\Loop}[1]{\def\args{#1}\ifx\args\empty\operatorname{\mathcal{L}}\else\operatorname{\mathcal{L}}(#1)\fi}
\newcommand{\TopB}{\Category{T}_{B}}

\newcommand{\catB}[1]{\operatorname{cat_{B}}(#1)}

\newcommand{\SuspB}[1]{\def\args{#1}\ifx\args\empty\operatorname{\mathcal{S}}_{B}\else\operatorname{\mathcal{S}_{B}}(#1)\fi}
\newcommand{\LoopB}[1]{\def\args{#1}\ifx\args\empty\operatorname{\mathcal{L}}_{B}\else\operatorname{\mathcal{L}_{B}}(#1)\fi}
\newcommand{\TopBb}{\Category{T}^{*}_{B}}

\newcommand{\catBb}[1]{\operatorname{cat^{*}_{B}}(#1)}

\newcommand{\SuspBb}[1]{\def\args{#1}\ifx\args\empty\operatorname{\mathcal{S}}^{*}_{B}\else\operatorname{\mathcal{S}^{*}_{B}}(#1)\fi}
\newcommand{\LoopBb}[1]{\def\args{#1}\ifx\args\empty\operatorname{\mathcal{L}}^{*}_{B}\else\operatorname{\mathcal{L}^{*}_{B}}(#1)\fi}
\newcommand{\TopBB}{\Category{T}^{B}_{B}}

\newcommand{\CatBB}[1]{\operatorname{Cat^{B}_{B}}(#1)}

\newcommand{\catBB}[1]{\operatorname{cat^{B}_{B}}(#1)}
\newcommand{\catlenBB}[1]{\operatorname{catlen^{B}_{B}}(#1)}
\newcommand{\wgtBB}[1]{\operatorname{wgt^{B}_{B}}(#1)}
\newcommand{\MwgtBB}[1]{\operatorname{Mwgt^{B}_{B}}(#1)}
\newcommand{\cuplenBB}[1]{\operatorname{cup^{B}_{B}}(#1)}
\newcommand{\SuspBB}[1]{\def\args{#1}\ifx\args\empty\operatorname{\mathcal{S}}^{B}_{B}\else\operatorname{\mathcal{S}}^{B}_{B}(#1)\fi}
\newcommand{\LoopBB}[1]{\def\args{#1}\ifx\args\empty\operatorname{\mathcal{L}}^{B}_{B}\else\operatorname{\mathcal{L}}^{B}_{B}(#1)\fi}
\title{%
Topological Complexity is a Fibrewise L-S Category%
}
\author[Iwase]{Norio Iwase}
\email{iwase@math.kyushu-u.ac.jp}
\author[Sakai]{Michihiro Sakai}
\email[Sakai]{sakai@kurume-nct.ac.jp}
\address[Iwase]{Faculty of Mathematics,
  Kyushu University,
  Fukuoka 810-8560, Japan}
\address[Sakai]{Kurume National College of Technology, 
  Fukuoka 830-8555, Japan.}
\begin{document}
%
%
\begin{abstract}
Topological complexity $\TC{B}$ of a space $B$ is introduced by M. Farber to measure how much complex the space is, which is first considered on a configuration space of a motion planning of a robot arm.
We also consider a stronger version $\TCM{B}$ of topological complexity with an additional condition: in a robot motion planning, a motion must be stasis if the initial and the terminal states are the same.
Our main goal is to show the equalities $\TC{B} = \catBb{\double{B}}+1$ and $\TCM{B} = \catBB{\double{B}}+1$, where $\double{B}=B{\times}B$ is a fibrewise pointed space over $B$ whose projection and section are given by $p_{\double{B}}=\proj_{2} : B{\times}B \to B$ the canonical projection to the second factor and $s_{\double{B}}=\Delta_{B} : B \to B{\times}B$ the diagonal.
In addition, our method in studying fibrewise L-S category is able to treat a fibrewise space with singular fibres.
\end{abstract}
%
%
%
%
%
%
\date{\today, [First draft]}
%
%
\maketitle
\section{Introduction}
%


We say a pair of spaces $(X,A)$ is an NDR pair or $A$ is an NDR subset of $X$, if the inclusion map is a (closed) cofibration, in other words, the inclusion map has the (strong) Str{\o}m structure (see page 22 in G. Whitehead \cite{Whitehead:elements}).
When the set of the base point of a space is an NDR subset, the space is called well-pointed.

Let us recall the definition of a sectional category (see James \cite{James:intro-fibrewise}) which is originally defined and called by Schwarz `genus'.
\begin{defn}[Schwarz \cite{Schwarz:genus}, James \cite{James:ls-category-survey}]
For a fibration $p : E \to X$, the sectional cateory $\secat{p}$ (= one less than the Schwarz genus $\genus{p}$) is the minimal number $m \geq 0$ such that there exists a cover of $X$ by $(m{+}1)$ open subsets $U_{i} \subset X$ each of which admits a continuous section $s_{i} : U_{i} \to E$.
\end{defn}
The topological complexity of a robot motion planning is first introduced by M. Farber \cite{Farber:motion-planning} in 2003 to measure the discontinuity of a robot motion planning algorithm searching also the way to minimise the discontinuity.
At a more general view point, Farber defined a numerical invariant $\TC{B}$ of any topological space $B$: let $\Path{B}$ be the space of all paths in $B$.
Then there is a Serre path fibration $\pi : \Path{B} \to B{\times}B$ given by $\pi(\ell)=(\ell(0),\ell(1))$ for $\ell \in \Path{B}$.
\begin{defn}[Farber]
For a space $B$, the topological complexity $\TC{B}$ is the minimal number $m \geq 1$ such that there exists a cover of $B{\times}B$ by $m$ open subsets $U_{i}$ each of which admits a continuous section $s_{i} : U_{i} \to \Path{B}$ for $\pi : \Path{B} \to B{\times}B$.
\end{defn}
%

By definition, we can observe that the topological complexity is nothing but the Schwartz genus or the sectional category.

%
Farber has further introduced a new invariant restricting motions by giving two additional conditions on the section $s : U \to \Path{B}$ (see Farber \cite{Farber:robot-motion-planning}).
%
\begin{enumerate}
\item
$s(b,b)=c_{b}$ the constant path at $b$ for any $b \in B$,
\item
$s(b_{1},b_{2})=s(b_{2},b_{1})^{-1}$ if $(b_{1},b_{2}) \in U$.
\end{enumerate}
It gives a stronger invariant than the topological complexity, and the $\integral/2$-equivariant theory must be applied as in Farber-Grant \cite{FG:smp}.
This new topological invariant, in turn, suggests us another motion planning under the condition that a motion is stasis if the initial and the terminal states are the same.
Let us state more precisely.
\begin{defn}
For a space $B$, the `monoidal' topological complexity $\TCM{B}$ is the minimal number $m \geq 1$ such that there exists a cover of $B{\times}B$ by $m$ open subsets $U_{i} \supset \Delta(B)$ each of which admits a continuous section $s_{i} : U_{i} \to \Path{B}$ for the Serre path fibration $\pi : \Path{B} \to B{\times}B$ satisfying $s_{i}(b,b)=c_{b}$ for any $b \in B$.
\end{defn}
\begin{rem}
This new topological complexity $\TCM{}$ is \textbf{not} a homotopy invariant, in general.
However, it is a homotopy invariant if we restrict our working category to the category of a space $B$ such that the pair $(B{\times}B,\Delta(B))$ is NDR.
\end{rem}

On the other hand, a fibrewise \textit{pointed} L-S category of a fibrewise pointed space is introduced and studied by James-Morris \cite{JM:fibrewise-category}.
Let us recall the definition:
\begin{defn}[James-Morris \cite{JM:fibrewise-category}]
\begin{enumerate}
\item
Let $X$ be a fibrewise pointed space over $B$.
The fibrewise \textbf{pointed} L-S category $\catBB{X}$ is the minimal number $m \geq 0$ such that there exists a cover of $X$ by $(m+1)$ open subsets $U_{i} \supset s_{X}(B)$ each of which is fibrewise null-homotopic in $X$ by a fibrewise pointed homotopy.
If there are no such $m$, we say $\catBB{X} = \infty$.
\item
Let $f : Y \to X$ be a fibrewise pointed map over $B$.
The fibrewise \textbf{pointed} L-S category $\catBB{f}$ is the minimal number $m \geq 0$ such that there exists a cover of $Y$ by $(m+1)$ open subsets $U_{i} \supset s_{Y}(B)$, where the restriction $f\vert_{U_{i}}$ to each subset is fibrewise compressible into $s_{X}(B)$ in $X$ by a fibrewise pointed homotopy.
If there are no such $m$, we say $\catBB{f} = \infty$.
\end{enumerate}
\end{defn}
To describe our main result, we further introduce a new unpointed version of fibrewise L-S category:
the fibrewise L-S category $\catB{ \ }$ of an fibrewise \textit{unpointed} space is also defined by James and Morris \cite{JM:fibrewise-category} as the minimum number (minus one) of open subsets which cover the given space and are fibrewise null-homotopic (see also James \cite{James:intro-fibrewise} and Crabb-James \cite{CJ:fibrewise}).
In this paper, we give a new version of a fibrewise \textit{unpointed} L-S category of a fibrewise \textit{pointed} space as follows:
\begin{defn}\label{def:fibrewise-unpointed-ls}
\begin{enumerate}
\item
Let $X$ be a fibrewise pointed space over $B$.
The fibrewise \textbf{unpointed} L-S category $\catBb{X}$ is the minimal number $m \geq 0$ such that there exists a cover of $X$ by $(m+1)$ open subsets $U_{i}$ each of which is fibrewise compressible into $s_{X}(B)$ in $X$ by a fibrewise homotopy.
If there are no such $m$, we say $\catBb{X} = \infty$.
\item
Let $f : Y \to X$ be a fibrewise pointed map over $B$.
The fibrewise \textbf{unpointed} L-S category $\catBb{f}$ is the minimal number $m \geq 0$ such that there exists a cover of $Y$ by $(m+1)$ open subsets $U_{i}$, where the restriction $f\vert_{U_{i}}$ to each subset is fibrewise compressible into $s_{X}(B)$ in $X$ by a fibrewise homotopy.
If there are no such $m$, we say $\catBb{f} = \infty$.
\end{enumerate}
\end{defn}

For a given space $B$, we define a fibrewise pointed space $\double{B}$ by $\double{B}=B{\times}B$ with $p_{\double{B}}=\proj_{2} : B{\times}B \to B$ and $s_{\double{B}}=\Delta_{B} : B \to B{\times}B$ the diagonal.
One of our main goals of this paper is to show the following theorem.
\begin{thm}\label{thm:main1}
For a space $B$, 
we have the following equalities.
\begin{enumerate}
\item
$\TC{B} = \catBb{\double{B}}+1$.
\item
$\TCM{B} = \catBB{\double{B}}+1$.
\end{enumerate}
\end{thm}

Farber and Grant has also introduced lower bounds for the topological complexity by using the cup length and category weight (see Rudyak \cite{Rudyak:ls-cat_mfds2,Rudyak:ls-cat_mfds3} for example) on the ideal of zero-divisors, i.e, the kernel of $\Delta^{*} : H^{*}(B{\times}B;R) \to H^{*}(B;R)$.

\begin{defn}[Farber \cite{Farber:motion-planning} and Farber-Grant  \cite{FG:smp}]\label{def:zero-cuplen}
For a space $B$ and a ring $R \ni 1$, the zero-divisors cup-length $\TCcuplen{B}{R}$ and the TC-weight $\TCwgt{u;R}{\pi}$ for $u \in I = \ker \Delta^{*} : H^{*}(B{\times}B;R) \to H^{*}(B;R)$ is defined as follows.
\begin{enumerate}
\item
$\TCcuplen{B}{R} = \Max\left\{m{\geq}0 \midvert H^{*}(B{\times}B;R) \supset I^{m}\not=0\right\}$
\item
$\TCwgt{u;R}{\pi} = \Max\left\{m{\geq}0 \midvert \forall{f : Y \to B{\times}B \ (\secat{f^{*}\pi} < m}), \ f^{*}(u)=0\right\}$
\end{enumerate}
\end{defn}
In the category $\TopBB$ of fibrewise pointed spaces with base space $B$ and maps between them, we also have corresponding definitions.
\begin{defn}
For a fibrewise pointed space $X$ over $B$ and a ring $R \ni 1$ and $u \in I = H^{*}(X,B;R) \subset H^{*}(X;R)$, we define
\begin{enumerate}
\item
$\cuplenBB{X;R} = \Max\left\{m{\geq}0 \midvert \exists\{u_{1},{\cdots},u_{m} \in H^{*}(X,B;R)\} \ \text{s.t.} \ u_{1}{\cdots}u_{m}\not=0\right\}$
\item
$\wgtBB{u;R} = \Max\left\{m{\geq}0 \midvert \forall{f : Y \to X \in \TopBB \ (\catBB{f} < m)}, \ f^{*}(u)=0\right\}$
\end{enumerate}
\end{defn}
This immediately implies the following.
\begin{thm}\label{thm:main2}
For a space $B$, we have $\TCcuplen{B}{R}=\cuplenBB{\double{B};R}$ for a ring $R \ni 1$.
\end{thm}
Motivating by this equality, we proceed to obtain the following result.
\begin{thm}\label{thm:main3}
For any space $B$,  
any element $u \in H^{*}(B{\times}B,\Delta(B);R)$ and a ring $R \ni 1$, 
we have $\TCwgt{u;R}{\pi} = \wgtBB{u;R}$.
\end{thm}
Let us consider one technical condition on a fibrewise pointed space:
\begin{thm}\label{thm:well-pointed}
For any space $B$ having the homotopy type of a locally finite simplicial complex, we may assume that $\double{B}$ is fibrewise well-pointed up to homotopy.
\end{thm}
The following is the main result of our paper.
\begin{thm}\label{thm:main4}
For any fibrewise \textit{well-pointed} space $X$ over $B$, we have $\catBB{X} = \catBb{X}$.
So, if $B$ is a locally finite simplicial complex, we have $\TC{B}=\TCM{B}=\catBB{\double{B}}+1$.
\end{thm}

In \cite{Sakai:thesis}, Sakai showed, in his study of the fibrewise \textit{pointed} L-S category of a fibrewise well-pointed spaces, using Whitehead style definition, that we can utilise $A_{\infty}$ methods used in the study of L-S category (see Iwase \cite{Iwase:counter-ls,Iwase:counter-ls-m}).
Let us state the Whitehead style definitions of fibrewise L-S categories following \cite{Sakai:thesis}.
\begin{defn}
Let $X$ be a fibrewise \textit{well-pointed} space over $B$.
The fibrewise \textbf{pointed} L-S category $\catBB{X}$ is the minimal number $m \geq 0$ such that the $(m{+}1)$-fold fibrewise diagonal $\Delta_{B}^{m+1} : X \to \overset{m+1}{\product_{B}}X$ is compressible into the fibrewise fat wedge $\overset{m+1}{\fatvee_{B}}X$ in $\TopBB$.
If there are no such $m$, we say $\catBB{X} = \infty$.
\end{defn}
We remark that this new definition coincides with the ordinary one, if the total space $X$ is a finite simplicial complex.

The above Whitehead-style definition allows us to define the module weight, cone length and categorical length, and moreover, to give their relationship as in Section \ref{sect:upper-lower}.
To show that, we need a criterion given by fibrewise $A_{\infty}$ structure on the fibrewise loop space (see Sections \ref{sect:a-infty}--\ref{sect:whitehead-defn}).

\section{Proof of Theorem \ref{thm:main1}}


First, we show the equality $\TCM{B}=\catBB{\double{B}}+1$: assume $\TCM{B}=m{+}1$, $m \geq 0$ and that there are an open cover $\bigcup^{m}_{i=0} U_{i} = B{\times}B$ and a series of sections $s_{i} : U_{i} \to \Path{B}$ of $\pi : \Path{B} \to \double{B}$ satisfying $s_{i}(b,b)=c_{b}$ for $b \in B$, since we are considering monoidal topological complexity.
Then each $U_{i}$ is fibrewise compressible relative to $\Delta(B)$ into $\Delta(B) \subset B{\times}B=\double{B}$ by a homotopy $H_{i} : U_{i}{\times}[0,1] \to B{\times}B$ given by the following:
$$
H_{i}(a,b;t)=(s_{i}(a,b)(t),b),\quad (a,b) \in U_{i}, \ t \in [0,1],
$$
where we can easily check that $H_{i}$ gives a fibrewise compression of $U_{i}$ relative to $\Delta(B)$ into $\Delta(B) \subset B{\times}B$.
Since $\bigcup_{i=0}U_{i}=B{\times}B=\double{B}$, we obtain $\catBB{\double{B}} \leq m$, and hence we have $\catBB{\double{B}}+1 \leq \TCM{B}$.

Conversely assume that $\catBB{\double{B}}=m$, $m \geq 0$ and there is an open cover $\bigcup^{m}_{i=0} U_{i} = \double{B}$ of $\double{B}=B{\times}B$ where $U_{i}$ is fibrewise compressible relative to $\Delta(B)$ into $\Delta(B) \subset \double{B}=B{\times}B$:
let us denote the compression homotopy of $U_{i}$ by $H_{i}(a,b;t) = (\sigma_{i}(a,b;t),b)$ for $(a,b) \in U_{i}$ and $t \in [0,1]$, where $\sigma_{i}(a,b;0)=a$ and $\sigma_{i}(a,b;1)=b$.
Hence we can define a section $s_{i} : U_{i} \to \Path{B}$ by the formula
$$
s_{i}(a,b)(t) = \sigma_{i}(a,b;t)\quad t \in [0,1].
$$
Since $\bigcup_{i=0}U_{i}=B{\times}B$, we obtain $\TCM{B} \leq m{+}1$ and hence we have $\TCM{B} \leq \catBB{\double{B}}+1$.
Thus we have $\TCM{B} = \catBB{\double{B}}+1$.

Second, we show the equality $\TC{B}=\catBb{\double{B}}+1$: assume $\TC{B}=m{+}1$, $m \geq 0$ and that there is a open cover $\bigcup^{m}_{i=0} U_{i} = B{\times}B$ and a section $s_{i} : U_{i} \to \Path{B}$ of $\pi : \Path{B} \to \double{B}$.
Then each $U_{i}$ is fibrewise compressible into $\Delta(B) \subset B{\times}B=\double{B}$ by a homotopy $H_{i} : U_{i}{\times}[0,1] \to B{\times}B$ which is given by 
$$
H_{i}(a,b;t)=(s(a,b)(t),b),\quad (a,b) \in U_{i}, \ t \in [0,1],
$$
where we can easily check that $H$ gives a fibrewise compression of $U_{i}$ into $\Delta(B) \subset B{\times}B=\double{B}$.
Since $\bigcup_{i=0}U_{i}=B{\times}B=\double{B}$, we obtain $\catBb{\double{B}} \leq m$, and hence we have $\catBb{\double{B}}+1 \leq \TC{B}$.

Conversely assume that $\catBb{\double{B}}=m$, $m \geq 0$ and there is an open cover $\bigcup^{m}_{i=0} U_{i} = \double{B}$ of $\double{B}=B{\times}B$ where $U_{i}$ is fibrewise compressible into $\Delta(B) \subset B{\times}B=\double{B}$:
the compression homotopy is described as $H_{i}(a,b;t) = (\sigma_{i}(a,b;t),b)$ for $(a,b) \in U_{i}$ and $t \in [0,1]$, such that $\sigma_{i}(a,b;0)=a$ and $\sigma_{i}(a,b;1)=b$.
Hence we can define a section $s_{i} : U_{i} \to \Path{B}$ by the formula
$$
s_{i}(a,b)(t) = \sigma_{i}(a,b;t)\quad t \in [0,1].
$$
Since $\bigcup_{i=0}U_{i}=B{\times}B$, we obtain $\TC{B} \leq m{+}1$ and hence we have $\TC{B} \leq \catBb{\double{B}}+1$.
Thus we have $\TC{B} = \catBb{\double{B}}+1$.
\qed

%
\section{Proof of Theorem \ref{thm:main3}}


Assume that $\wgtBB{u;R} = m$, where $u \in H^{*}(B{\times}B,\Delta(B))$ and $f : Y \to \double{B}=B{\times}B$ a map of $\secat{f^{*}\pi} < m$.
Then there is an open cover $\bigcup^{m}_{i=1} U_{i} = Y$ and a series of maps $\{\sigma_{i} : U_{i} \to \Path{B}\,;\, 1 \leq i \leq m\}$ satisfying $\pi{\comp}\sigma_{i}=f\vert_{U_{i}}$.
Let $\hat{Y} = Y \amalg B$ with projection $p_{\hat{Y}}$ and section $s_{\hat{Y}}$ given by
\begin{align*}&
p_{\hat{Y}}\vert_{Y} = p_{Y}, \quad p_{\hat{Y}}\vert_{B} = \id_{B}\quad\text{and}\quad
s_{\hat{Y}} : B \hookrightarrow Y \amalg B = \hat{Y}.
\end{align*}
Then we can extend $f$ to a map $\hat{f} : \hat{Y} \to \double{B}$ by the formula 
$$
\hat{f}\vert_{Y} = f, \quad \hat{f}\vert_{B} = s_{\double{B}} = \Delta.
$$
By putting $\hat{U}_{i}=U_{i} \amalg B$ which is open in $\hat{Y}$, we obtain an open cover $\bigcup^{m}_{i=1} \hat{U}_{i} = \hat{Y}$ and a series of maps $\hat{\sigma}_{i} : \hat{U}_{i} \to \Path{B}$ satisfying $\pi{\comp}\hat{\sigma}_{i}=\hat{f}\vert_{\hat{U}_{i}}$:
$$
\hat{\sigma}_{i}\vert_{U_{i}} = \sigma_{i},\quad \hat{\sigma}_{i}\vert_{B} = s_{\Path{B}}.
$$
Hence there is a fibrewise homotopy $\Phi_{i} : \hat{U}_{i}{\times}[0,1] \to \double{B}$ such that $\Phi_{i}(y,0)=\hat{f}(y)$ and $\Phi_{i}(y,1) \in \Delta(B)$ given by the following formula.
$$
\Phi_{i}(y,t) = (\hat{\sigma}_{i}(y)(t),\hat{\sigma}_{i}(y)(1)), \quad (y,t) \in \hat{U}_{i}{\times}[0,1],
$$
so that we have $\Phi_{i}(y,0)=(\hat{\sigma}_{i}(y)(0),\hat{\sigma}_{i}(y)(1))=\pi{\comp}\hat{\sigma}_{i}(y)=\hat{f}(y)$ and $\Phi_{i}(y,1)=(\hat{\sigma}_{i}(y)(1),\hat{\sigma}_{i}(y)(1)) \in \Delta(B)$.
Moreover, for any $(b,t) \in B{\times}[0,1]$, we have $\Phi_{i}(b,t)=(\hat{\sigma}_{i}(b)(t),\hat{\sigma}_{i}(b)(1))=(s_{\Path{B}}(t),s_{\Path{B}}(1))=(b,b)$.
Thus $\Phi_{i}$ gives a fibrewise pointed compression homotopy of $\hat{f}\vert_{\hat{U}_{i}}$ into $\Delta(B)$.
Then it follows that $\catBB{\hat{f}} < m$ and hence we obtain $f^{*}(u)=0$ and $\TCwgt{u;R}{\pi} \geq m$.
Thus we obtain $\TCwgt{u;R}{\pi} \geq m=\wgtBB{u;R}$.

Conversely assume that $\TCwgt{u;R}{\pi} = m$, where $u \in H^{*}(B{\times}B,\Delta(B))$ and $f : Y \to B{\times}B$ such that $\catBB{f} < m$.
Then there exists an open covering $\bigcup^{m}_{i=1} U_{i} = Y$ with $U_{i} \supset s_{Y}(B)$ and a sequence of fibrewise homotopies $\left\{\phi_{i} : U_{i}{\times}[0,1] \to B{\times}B \right\}$ such that $\phi_{i}(y,0) = f\vert_{U_{i}}(y)$, $\phi_{i}(y,1) \in \Delta(B)$ and $\proj_{2}{\comp}\phi_{i}(y,t)=\proj_{2}{\comp}f(y)$ for $(y,t) \in U_{i}{\times}[0,1]$.
Hence there is a sequence of maps $\{\sigma_{i} : U_{i} \to \Path{B}\}$ given by 
$$
\sigma_{i}(y)(t) = \proj_{1}{\comp}\phi_{i}(y,t), \quad y \in U_{i}, \ t \in [0,1]
$$
such that $\pi{\comp}\sigma_{i}(y)=(\proj_{1}{\comp}\phi_{i}(y,0),\proj_{1}{\comp}\phi_{i}(y,1))=f(y)$ since $\proj_{2}{\comp}\phi_{i}(y,t)=\proj_{2}{\comp}f(y)$ for $(y,t) \in U_{i}{\times}[0,1]$.
Thus we obtain $\secat{f^{*}\pi} < m$, and hence $f^{*}(u) = 0$.
This implies $\wgtBB{u;R} \geq m=\TCwgt{u;R}{\pi}$ and hence $\wgtBB{u;R} = \TCwgt{u;R}{\pi}$.
\qed

%
\section{Proof of Theorem \ref{thm:well-pointed}}
%
The proof of Lemma 2 in \S2 of Milnor \cite{Milnor:homotopy-type} implies the following:
\begin{lem}\label{lem:milnor}
The pair $(B{\times}B,\Delta(B))$ is an NDR-pair.
\end{lem}
\begin{proof}
For each vertex $\beta$ of $B$, let $V_{\beta}$ be the star neighbourhood in $B$ and $V = \bigcup_{\beta} V_{\beta}{\times}V_{\beta} \subset B{\times}B$.
Then the closure $\bar{V} = \bigcup_{\beta} \bar{V}_{\beta}{\times}\bar{V}_{\beta}$ is a subcomplex of $B{\times}B$.
For the barycentric coordinates $\{\xi_{\beta}\}$ and $\{\eta_{\beta}\}$ of $x$ and $y$, resp, we see that $(x,y) \in V$ if and only if $\sum_{\beta}\Min(\xi_{\beta},\eta_{\beta})>0$ and that $\sum_{\beta}\Min(\xi_{\beta},\eta_{\beta})=1$ if and only if the barycentric coordinates of $x$ and $y$ are the same, or equivalently, $(x,y) \in \Delta(B)$.
Hence we can define a continuous map $v : B{\times}B \to [-1,1]$ by the following formula.
\begin{align*}
v(x,y)=\begin{cases}
2\sum_{\beta}\Min(\xi_{\beta},\eta_{\beta})-1, & \text{if $(x,y) \in \bar{V}$,}\\
-1, & \text{if $(x,y) \not\in V$.}
\end{cases}
\end{align*}
Then we have that $v^{-1}(1) = \Delta(B)$.
Let $U=v^{-1}((0,1])$ an open neighbourhood of $\Delta(B)$.
Using Milnor's map $s$, we obtain a pair of maps $(u,h)$ as follows:
\begin{align*}&
u(x,y) = \Min\{1,1{-}v(x,y)\}\quad\text{and}
\\&
h(x,y,t) = 
(s(x,y)(\Min\{t,w(x,y)\}),y), 
\end{align*}
where $w(x,y)=u(x,y)+v(x,y)=\Min\{1,1{+}v(x,y)\}$.
Note that $w(x,y)=1$ if $(x,y) \in U$ and that $w(x,y)=0$ if $(x,y) \not\in V$.
Then $u^{-1}(0)=\Delta(B)$, $u^{-1}([0,1))=U$ and $h(x,y,1) = (y,y) \in \Delta(B)$ if $(x,y) \in U$.
Moreover, $\proj_{2}{\comp}h(x,y,t) = y$ and $h(x,x,t) = (s(x,x)(t),x)=(x,x)$ for any $x,y \in B$ and $t \in [0,1]$.
Thus the data $(u,h)$ gives the fibrewise Str{\o}m structure on $(B{\times}B,\Delta(B))$.
\end{proof}

%
\section{Proof of Theorem \ref{thm:main4}}

Let $X$ be a fibrewise well-pointed space over $B$ and $\hat{X}$ the fiberwise pointed space obtained from $X$ by giving a fibrewise whisker.
More precisely, we define $\hat{X}$ be the mapping cylinder of $s_{X}$, 
$$
\hat{X} = X \cup_{s_{X}} B{\times}[0,1],\quad X \ni s_{X}(b) \sim (b,0) \in B{\times}[0,1] \ \text{for any $b \in B$,}
$$
with projection $p_{\hat{X}}$ and section $s_{\hat{X}}$ given by the formulas
\begin{align*}&
p_{\hat{X}}\vert_{X}=p_{X},\quad p_{\hat{X}}\vert_{B{\times}[0,1]}(b,t)=b, \quad \text{for $(b,t) \in B{\times}[0,1]$},
\\&
s_{\hat{X}}(b)=(b,1) \in B{\times}[0,1] \subset \hat{X}.
\end{align*}
Then by the definition of Str{\o}m structure, $X$ is fibrewise pointed homotopy equivalent to $\hat{X}$ the fibrewise whiskered space over $B$.
So we have $\catBB{X}=\catBB{\hat{X}}$ and $\catBb{X}=\catBb{\hat{X}}$.

Assume that $\catBB{X}=m \geq 0$.
Then it is clear by definition that $\catBb{X} \leq m=\catBB{X}$.

Conversely assume that $\catBb{X}=m \geq 0$.
Then there is an open cover $\bigcup^{m}_{i=0}U_{i} = X$ such that $U_{i}$ is compressible into $s_{X}(B) \subset X$.
Hence there is a fibrewise homotopy $\Phi_{i} : U_{i}{\times}[0,1] \to X$ such that $\Phi_{i}(x,0)=x$, $\Phi_{i}(x,1)=s_{X}(p_{X}(x))$ and $p_{X}{\comp}\Phi_{i}(x,t)=p_{X}(x)$.
We define $\hat{U}_{i}$ as follows:
$$
\hat{U}_{i} = U_{i} \cup_{s_{X}} (s_{X})^{-1}(U_{i}){\times}[0,1] \cup B{\times}(\frac{2}{3},1].
$$
We also define a fibrewise pointed homotopy $\hat{\Phi}_{i} : \hat{U}_{i}{\times}[0,1] \to \hat{X}$ as follows:
$$
\hat{\Phi}_{i}(\hat{x},t) = \begin{cases}
\Phi_{i}(x,t), &\hat{x}=x \in X,
\\
\Phi_{i}(s_{X}(b),t{-}3s), &\hat{x}=(b,s) \in (s_{X})^{-1}(U_{i}){\times}(0,\frac{t}{3}),
\\
s_{X}(b), &\hat{x}=(b,\frac{t}{3}), b \in (s_{X})^{-1}(U_{i}),
\\
(b,\frac{6s{-}2t}{6{-}3t}), &\hat{x}=(b,s) \in (s_{X})^{-1}(U_{i}){\times}(\frac{t}{3},\frac{2}{3}),
\\
(b,\frac{2}{3}), &\hat{x}=(b,\frac{2}{3}), b \in (s_{X})^{-1}(U_{i}),
\\
(b,s), &\hat{x}=(b,s) \in B{\times}(\frac{2}{3},1].
\end{cases}
$$
It is then easy to see that $\hat{U}_{i}$'s cover the entire $X$, and hence we have $\catBB{\hat{X}} \leq m=\catBb{X}$.
Thus $\catBB{X} \leq \catBb{X}$ and hence $\catBB{X}=\catBb{X}$.
In particular, we have $\TC{B}=\TCM{B}$ for a locally finite simplicial complex $B$.
\qed

\section{Fibrewise $A_{\infty}$ structures}\label{sect:a-infty}

From now on, we work in the category $\TopBB$.
For any $X$ a fibrewise pointed space over $B$, we denote by $p_{X} : X \to B$ its projection and by $s_{X} : B \to X$ its section.

We say that a pair $(X,A)$ of fibrewise pointed spaces over $B$ is a fibrewise NDR-pair or that $A$ is a fibrewise NDR subset of $X$, if the inclusion map $A \hookrightarrow X$ is a fibrewise cofibration, in other words, the inclusion has the fibrewise (strong) Str{\o}m structure (see Crabb-James \cite{CJ:fibrewise}).
Since $B$ is the zero object in $\TopBB$, for any given fibrewise pointed space $X$ over $B$, we always have a pair $(X,B)$ in $\TopBB$, where we regard $s_{X}(B)=B$.
When the pair $(X,B)$ is fibrewise NDR, the space $X$ is called fibrewise well-pointed.

\begin{prop}[Crabb-James \cite{CJ:fibrewise}]
\begin{enumerate}
\item
If $(X,A)$ and $(X',A')$ are fibrewise NDR-pairs, then so is $(X,A){{\times}_{B}}(X',A')=(X{{\times}_{B}}X',X{{\times}_{B}}A'{\cup}A{{\times}_{B}}X')$.
\item
If $(X,A)$ is a fibrewise NDR-pair, then so is $(\overset{m}{\product}_{B}X,\overset{m}{\fatvee}_{B}(X,A))$, which is defined by induction for all $m \geq 1$:
\begin{align*}&
(\overset{1}{\product}_{B}X,\overset{1}{\fatvee}_{B}(X,A)) = (X,A),
\\&
(\overset{m+1}{\product_{B}}X,\overset{m+1}{\fatvee_{B}}(X,A)) = (\overset{m}{\product}_{B}X,\overset{m}{\fatvee}_{B}(X,A)){\times}_{B}(X,A).
\end{align*}
\end{enumerate}
\end{prop}
If $X$ is a fibrewise pointed space over $B$, then by taking $A=B$, we obtain a fibrewise subspace $\overset{m+1}{\fatvee_{B}}(X,B)$ of $\overset{m+1}{\fatvee_{B}}X$, which is called an $(m{+}1)$-fold fibrewise fat-wedge of $X$, and is often denoted by $\overset{m+1}{\fatvee_{B}}X$.
In addition, the pair $(\overset{m+1}{\product_{B}}X,\overset{m+1}{\fatvee_{B}}X)$ is a fibrewise NDR-pair for all $m \geq 0$, if $X$ is fibrewise well-pointed.

\begin{expls}
\begin{enumerate}
\item\label{ex:trivial}
Let $X$ be a fibrewise pointed space over $B$ with $p_{X}=\proj_{2} : X=F{\times}B \to B$ the canonical projection to the second factor and $s_{X}=\incl_{2} : B \hookrightarrow F{\times}B=E$ the canonical inclusion to the second factor.
Then $X$ is a fibrewise pointed space over $B$.
\item\label{ex:diagonal}
Let $X=B{\times}B$, $p_{X}=\proj_{2} : B{\times}B \to B$ the canonical projection to the second factor and $s_{X}=\Delta_{B} : B \hookrightarrow B{\times}B$ the diagonal.
Then $X$ is a fibrewise pointed space over $B$.
\item\label{ex:group}
Let $G$ be a topological group, $EG$ the infinite join of $G$ with right $G$ action and $BG=EG/G$ the classifying space of $G$.
By considering $G$ as a left $G$ space by the adjoint action, we obtain a fibrewise pointed space $X=EG \times_{G} G$ with $p_{X} : EG \times_{G} G \to BG$ with section $s_{X} : BG \hookrightarrow EG \times_{G} \{e\} \subseteq EG \times_{G} G$.
\item\label{ex:free-loop}
Let $B$ be a space, $X=\fLoop{B}$ the space of free loops on $B$.
Then $p_{X} : \fLoop{B} \to B$ the evaluation map at $1 \in S^{1} \subset \complex$ is a fibration with section $s_{X} : B \to \fLoop{B}$ given by the inclusion of constant loops.
In view of Milnor's arguments, this example is homotopically equivalent to the example (\ref{ex:group}).
\end{enumerate}
\end{expls}

\begin{defn}
Let $\PathB{X} = \left\{\ell : [0,1] \to X \midvert \exists_{b \in B} \ \text{s.t.} \ \forall_{t \in [0,1]} \ p_{X}(\ell(t)){=}b\right\}$ the fibrewise free path space, $\LoopB{X} = \left\{\ell \in \PathB{X} \midvert \ell(1){=}\ell(0)\right\}$ the fibrewise free loop space and $\LoopBB{X} = \left\{\ell \in \PathB{X} \midvert \ell(1){=}\ell(0){=}s_{X}{\comp}p_{X}(\ell(0))\right\}$ the fibrewise pointed loop space.
For any $m \geq 0$, we define an $A_{\infty}$ structure of $\LoopBB{X}$ as follows.
\begin{enumerate}
\item
$E_{B}^{m+1}(\LoopBB{X})$ as the homotopy pull-back in $\TopBB$ of 
$B \hookrightarrow \overset{m+1}{\product_{B}}X \hookleftarrow \overset{m+1}{\fatvee_{B}}X$,
\item
$P_{B}^{m}(\LoopBB{X})$ as the homotopy pull-back in $\TopBB$ of $X \xrightarrow{\Delta_{B}^{m+1}} \overset{m+1}{\product_{B}}X \hookleftarrow \overset{m+1}{\fatvee_{B}}X$,
\item
$e^{X}_{m} : P_{B}^{m}(\LoopBB{X}) \to X$ as the induced map from the inclusion $\overset{m+1}{\fatvee_{B}}X \hookrightarrow \overset{m+1}{\product_{B}}X$ by the diagonal $\Delta_{B}^{m+1} : X \to \overset{m+1}{\product_{B}}X$ and 
\item
$p_{B}^{\LoopBB{X}} : E_{B}^{m+1}(\LoopBB{X}) \to P_{B}^{m}(\LoopBB{X})$ as a map of fibrewise pointed spaces induced from the section $s_{X} : B \to X$, since the section $B \hookrightarrow \overset{m+1}{\product_{B}}X$ is nothing but the composition $\Delta_{B}^{m+1}{\comp}s_{X} : B \xrightarrow{s} X \xrightarrow{\Delta_{B}^{m+1}} \overset{m+1}{\product_{B}}X$.
\end{enumerate}
\end{defn}

We further investigate to understand an $A_{\infty}$ stucture in a fiberwise view point, using fibrewise constructions.
Clearly, these constructions are \textit{not} exactly the Ganea-type fibre-cofibre constructions but the following.
\begin{prop}[Sakai]\label{prop:pushpull}
Let $X$ be a fibrewise pointed space over $B$ and $m \geq 0$.
Then $P_{B}^{m+1}(\LoopBB{X})$ has the homotopy type of a push-out of $p_{B}^{\LoopBB{X}} : E_{B}^{m+1}(\LoopBB{X})$ $\to$ $P_{B}^{m}(\LoopBB{X})$ and the projection $E_{B}^{m+1}(\LoopBB{X}) \to B$.
\end{prop}
This is a direct consequence of the following lemma.
\begin{lem}\label{lem:pushpull}
Let $(X,A)$ and $(X',A')$ be fibrewise NDR-pairs of fibrewise pointed spaces over $B$ and $Z$ a fibrewise pointed space over $B$ with fibrewise maps $f : Z \to X$ and $g : Z \to X'$.
Then the homotopy pull-back $\Omega_{(f,g),k}$ of maps $(f,g) : Z \to X{{\times}_{B}}X'$ and $k : X{\times}_{B}A' \cup A{\times}_{B}X' \hookrightarrow X{\times}_{B}X'$ has naturally the homotopy type of the reduced homotopy push-out $W = \Omega_{g,j} \cup_{p_{2}} \left\{\Omega_{(f,g),i{\times}j}\wedge_{B}(B{\times}J^{+})\right\} \cup_{p_{1}} \Omega_{f,i}$ of $p_{1} : \Omega_{(f,g),i{\times}j} \to \Omega_{f,i}$ and $p_{2} : \Omega_{(f,g),i{\times}j} \to \Omega_{g,j}$, where $J=[-1,1]$ and
\begin{align*}&
\Omega_{(f,g),k} = \left\{(z,\ell,\ell') \in Z {\times_{B}} \PathB{X} {\times_{B}} \PathB{X'} \midvert \substack{
f(z)=\ell(0), \ g(z)=\ell'(0),\\ (\ell(1),\ell'(1)) \in A{\times}_{B}X' \cup X{\times}_{B}A'}\right\},
\\&
\Omega_{(f,g),i{\times}j} = \left\{(z,\ell,\ell') \in \Omega_{(f,g),k} \midvert (\ell(1),\ell'(1)) \in A{\times}_{B}A'\right\},
\\&
\Omega_{f,i} = \left\{(z,\ell) \in Z {\times_{B}} \PathB{X} \midvert f(z){=}\ell(0), \ \ell(1) {\in} A\right\},
\\&
\Omega_{g,j} = \left\{(z,\ell') \in Z {\times_{B}} \PathB{X'} \midvert g(z){=}\ell'(0), \ \ell'(1) {\in} A'\right\},
\end{align*}
$p_{1}(z,\ell,\ell')=(z,\ell)$ and $p_{2}(z,\ell,\ell')=(z,\ell')$.
\end{lem}
\begin{proof}[Outline of the proof]
The proof of Lemma \ref{lem:pushpull} is quite similar to that of Theorem 1.1 in Sakai \cite{Sakai:push_pull} (which is based on Iwase \cite{Iwase:counter-ls}) by replacing $(Y,B)$ in \cite{Sakai:push_pull} by $(X',A')$, defining and using the following spaces.
\begin{align*}&
\widehat{W} = \Omega_{(f,g),i{\times}\id_{X'}}{\times}\{-1\} \cup \left\{\Omega_{(f,g),i{\times}j}{\times}J\right\} \cup \Omega_{(f,g),\id_{X}{\times}j}{\times}\{1\} 
\subset \Omega_{(f,g),k}{\times}J,
\\&
\Omega_{(f,g),\id_{X}{\times}j} = \left\{(z,\ell,\ell') \in \Omega_{(f,g),k} \midvert (\ell(1),\ell'(1)) \in X{\times}_{B}A'\right\},
\\&
\Omega_{(f,g),i{\times}\id_{X'}} = \left\{(z,\ell,\ell') \in \Omega_{(f,g),k} \midvert (\ell(1),\ell'(1)) \in A{\times}_{B}X'\right\}.
\end{align*}
The precise construction of homotopy equivalences and homotopies is identical to that in \cite{Sakai:push_pull} and is left to the readers.
\end{proof}

\begin{thm}
Let $X$ be a fibrewise well-pointed space over $B$.
Then the sequence $\{p_{B}^{\LoopBB{X}} : E_{B}^{m+1}(\LoopBB{X}) \to P_{B}^{m}(\LoopBB{X})\}$ gives a fibrewise pointed version of $A_{\infty}$-structure on the fibrewise pointed loop space $\LoopBB{X}$.
\end{thm}
Thus in the case when $X$ is a fibrewise well-pointed space over $B$, we assume that $P_{B}^{m}(\LoopBB{X})$ is an increasing sequence given by homotopy push-outs with a fibrewise fibration $e_{m}^{X} : P_{B}^{m}(\LoopBB{X}) \to X$ such that $e_{1}^{X} : \SuspBB{\LoopBB{X}} \to X$ is a fibrewise evaluation.
\begin{expls}
\begin{enumerate}
\item
Let $X$ be a fibrewise pointed space over $B$ with $p_{X}=\proj_{2} : F{\times}B \to B$ the canonical projection and $s_{X}=\incl_{2} : B \hookrightarrow F{\times}B$ the canonical inclusion.
Then $\LoopBB{X}=\Loop{F}{\times}B$ is given by $p_{\LoopBB{X}}=\proj_{2} : \Loop{F}{\times}B \to B$ and $s_{\LoopBB{X}}=\incl_{2} : B \hookrightarrow \Loop{F}{\times}B$.
\item
Let $X=B{\times}B$ be a fibrewise pointed space over $B$ with $p_{X}=\proj_{2} : B{\times}B \to B$ and $s_{X}=\Delta_{B} : B \hookrightarrow B{\times}B$ the diagonal.
Then $\LoopBB{X}=\fLoop{B}$ the free loop space on $B$, $p_{\LoopBB{X}} : \fLoop{B} \to B$ the evalation map at $1 \in S^{1} \subset \complex$ and $s_{\LoopBB{X}} : B \hookrightarrow \fLoop{B}$ the inclusion of constant loops.
\end{enumerate}
\end{expls}
\begin{rem}
When $E$ is a cell-wise trivial fibration on a polyhedron $B$ (see \cite{IS:func_qfib}), we can see that the canonical map $e_{\infty}^{E} : P_{B}^{\infty}(\LoopBB{E}) \to E$ is a homotopy equivalence by a similar arguments given in the proof of Theorem 2.9 of \cite{IS:func_qfib}.
\end{rem}

%
\section{Fibrewise L-S categories of fibrewise pointed spaces}\label{sect:whitehead-defn}

The fibrewise \textit{pointed} L-S category of an fibrewise pointed space is first defined by James and Morris \cite{JM:fibrewise-category} as the least number (minus one) of open subsets which cover the given space and are contractible by a homotopy fixing the base point in each fibre (see also James \cite{James:intro-fibrewise} and Crabb-James \cite{CJ:fibrewise}) and is redefined by Sakai in \cite{Sakai:thesis} as follows:
let $X$ be a fibrewise pointed space over $B$.
For given $k \geq 0$, we denote by $\overset{k+1}{\product_{B}}X$ the $(k{+}1)$-fold fibrewise product and by $\overset{k+1}{\fatvee_{B}}X$ the $(k{+}1)$-fold fibrewise fat wedge.
Then $\catBB{X} \leq m$ if the $(m{+}1)$-fold fibrewise diagonal map $\Delta_{B}^{m+1} : X \to \overset{m+1}{\product_{B}}X$ is compressible into the fibrewise fat wedge $\overset{m+1}{\fatvee_{B}}X$ in $\TopBB$.
If there is no such $m$, we say $\catBB{X} = \infty$.
Let us consider the case when $\catBB{X} < \infty$.
The definition of a fibrewise $A_{\infty}$ structure yields the following criterion.

\begin{thm}
Let $X$ be a fibrewise pointed space over $B$ and $m \geq 0$.
Then $\catBB{X} \leq m$ if and only if $\id_{X} : X \to X$ has a lift to $P_{B}^{m}(\LoopBB{X}) \overset{e_{m}^{X}}{\rightarrow} X$ in $\TopBB$.
\end{thm}
\begin{proof}
If $\catBB{X} \leq m$, then the fibrewise diagonal $\Delta^{m+1}_{B} : X \to \overset{m+1}{\product_{B}}X$ is compressible into the fibrewise fat wedge $\overset{m+1}{\fatvee_{B}}X \subset \overset{m+1}{\product_{B}}X$ in $\TopBB$.
Hence there is a map $\sigma : X \to P_{B}^{m}(\LoopBB{X})$ in $\TopBB$ such that $e^{X}_{m}{\comp}\sigma \sim_{B} 1_{X}$ in $\TopBB$.
The converse is clear by the definition of $P_{B}^{m}(\LoopBB{X})$.
\end{proof}

%
%

In the rest of this section, we work within the category $\TopB$ of fibrewise \textit{unpointed} spaces and maps between them.
But we concentrate ourselves to consider its full subcategory $\TopBb$ of all fibrewise pointed spaces, so in $\TopBb$, we have more maps than in $\TopBB$ while we have just the same objects as in $\TopBB$.

Let $X$ be a fibrewise pointed space over $B$.
For given $k \geq 0$, we denote by $\overset{k+1}{\product_{B}}X$ the $(k{+}1)$-fold fibrewise product and by $\overset{k+1}{\fatvee_{B}}X$ the $(k{+}1)$-fold fibrewise fat wedge.
Then $\catBb{X} \leq m$ if the $(m{+}1)$-fold fibrewise diagonal map $\Delta_{B}^{m+1} : X \to \overset{m+1}{\product_{B}}X$ is compressible into the fibrewise fat wedge $\overset{m+1}{\fatvee_{B}}X$ in $\TopBb$.
If there is no such $m$, we say $\catBb{X} = \infty$.
Let us consider the case when $\catBb{X} < \infty$.
The definition of a fibrewise $A_{\infty}$ structure yields the following.
\begin{thm}
Let $X$ be a fibrewise pointed space over $B$ and $m \geq 0$.
Then $\catBb{X} \leq m$ if and only if $\id_{X} : X \to X$ has a lift to $P_{B}^{m}(\LoopBB{X}) \overset{e_{m}^{X}}{\rightarrow} X$ in the category $\TopBb$.
\end{thm}
\begin{proof}
If $\catBb{X} \leq m$, then the fibrewise diagonal $\Delta^{m+1}_{B} : X \to \overset{m+1}{\product_{B}}X$ is compressible into the fibrewise fat wedge $\overset{m+1}{\fatvee_{B}}X \subset \overset{m+1}{\product_{B}}X$ in $\TopBb$.
Hence there is a map $\sigma : X \to P_{B}^{m}(\LoopBB{X})$ in $\TopBb$ such that $e^{X}_{m}{\comp}\sigma \sim_{B} 1_{X}$ in $\TopBb$.
The converse is clear by the definition of $P_{B}^{m}(\LoopBB{X})$.
\end{proof}

\section{Upper and lower estimates}\label{sect:upper-lower}

For $X$ a fibrewise pointed space over $B$, we define a fibrewise version of Ganea's strong L-S category (see Ganea \cite{Ganea:category}) of $X$ as $\CatBB{X}$ 
and also a fibrewise version of Fox's categorical length (see Fox \cite{Fox:ls-cat} and Iwase \cite{Iwase:cat-length}) of $X$ as $\catlenBB{X}$. 
\begin{defn}
Let $X$ be a fibrewise pointed space over $B$.
\begin{enumerate}
\item\label{def:strong-catBB}
$\CatBB{X}$ is the least number $m \geq 0$ such that there exists a sequence $\left\{(X_{i},h_{i})\midvert h_{i} : A_{i} {\to} X_{i-1}, \ 0{\leq}i{\leq}m\right\}$ of pairs of space and map satisfying $X_{0} = B$ and $X_{m} \simeq_{B} X$ in $\TopBB$ with the following homotopy push-out diagrams:
\begin{equation*}
\begin{diagram}
\node{A_{i}}
\arrow{e,t}{p_{A_{i}}}
\arrow{s,l}{h_{i}}
\node{B}
\arrow{s,r,J}{s_{X_{i}}}
\\
\node{X_{i-1}}
\arrow{e,J}
\node{X_{i}}
\end{diagram}
\end{equation*}
\item\label{def:cat-lengthBB}
$\catlenBB{X}$ is the least number $m \geq 0$ such that there exists a sequence $\left\{X_{i} \midvert h_{i} : A_{i} {\to} X_{i-1}, \ 0 {\leq} i {\leq} m\right\}$ of spaces satisfying $X_{0} = B$ and $X_{m} \simeq_{B} X$ in $\TopBB$ and that $\Delta_{B} : X_{i} \to X_{i}{\times}_{B}X_{i}$ is compressible into $X_{i}{\times}_{B}X_{i-1} \cup B{\times}_{B}X_{i}$ in $X_{m}{\times}_{B}X_{m}$.
\end{enumerate}
\end{defn}

A lower bound for the fibrewise L-S category of a fibrewise pointed space $X$ over $B$ can be described by a variant of cup length:
since $X$ is a fibrewise pointed space over $B$, there is a projection $p_{X} : X \to B$ with its section $s_{X} : B \to X$.
Hence we can easily observe for any multiplicative cohomology theory $h$ that 
\begin{align*}&
h^{*}(X) \cong h^{*}(B){\oplus}h^{*}(X,B),
\end{align*}
where we may identify $h^{*}(X,B)$ with the ideal $\ker s_{X}^{*} : h^{*}(X) \to h^{*}(B)$.
%
\begin{defn}
For a fibrewise pointed space $X$ over $B$ and any multiplicative cohomology theory $h$, we define
\begin{align*}&
\cuplenBB{X;h} = \Max\left\{m {\geq} 0 \midvert \exists\{u_{1},{\cdots},u_{m} \in h^{*}(X,B)\} \ \text{s.t.} \ u_{1}{\cdots}u_{m}\not=0\right\},
\\&
\cuplenBB{X} = \Max\left\{
\cuplenBB{X;h} \midvert \text{$h$ is a multiplicative cohomology theory}
\right\}.
\end{align*}
\end{defn}
We often denote $\cuplenBB{ \ ;h}$ by $\cuplenBB{ \ ;R}$ when $h^{\ast}( \ ) = H^{\ast}( \ ;R)$, where $R$ is a ring with unit.
%

Let us recall that the relationship between an $A_{\infty}$-structure and a Lusternik-Schnirelmann category gives the key observation in \cite{Iwase:counter-ls,Iwase:counter-ls-m,Iwase:ls-cat-bundle}.

On the other hand, Rudyak \cite{Rudyak:ls-cat_mfds2} and Strom \cite{Strom:essential-cat-wgt} introduced a homotopy theoretical version of Fadell-Husseini's category weight, which can be translated into our setting as follows: for any fibrewise pointed space $X$ over $B$, let $\{p^{\LoopBB{X}}_k \,{:}\, E_{B}^{k}({\LoopBB{X}})$ $\to$ $P_{B}^{k-1}(\LoopBB{X}) \,;\, k {\geq} 1\}$ be the fibrewise $A_{\infty}$-structure of $\LoopBB{X}$ in the sense of Stasheff \cite{Stasheff:higher-associativity} (see also \cite{IM:higher-associativity} for some more properties).
Let $h$ be a generalisd cohomology theory.
\begin{defn}\label{def:weightB}
For any $u \in h^{\ast}(X,B)$, we define
$$\wgtBB{u;h} = \Min\left\{m {\geq} 0%
\left\vert\,%
(e^{X}_{m})^{\ast}(u) \not= 0
\right.\right\},
$$
where $e^{X}_{m}$ is the composition of fibrewise maps $P_{B}^{m}(\LoopBB{X}) \hookrightarrow P_{B}^{\infty}(\LoopBB{X}) \underset{\simeq_{B}}{\xrightarrow{e^{X}_{\infty}}} X$.
\end{defn}
Using this, we introduce some more invariants as follows.
\begin{defn}
For any fibrewise pointed space $X$ over $B$, we define 
\begin{align*}&
\TCwgt{X;h}{\pi} = \Max\left\{
\TCwgt{u;h}{\pi}\,\vert\,u \in h^{*}(X,B)
\right\},
\\&
\TCwgt{X}{\pi} = \Max\left\{
\TCwgt{X;h}{\pi} \midvert \text{$h$ is a generalised cohomology theory}
\right\},
\\&
\wgtBB{X;h} = \Max\left\{
\wgtBB{u;h}\,\vert\,u \in h^{*}(X,B)
\right\},
\\&
\wgtBB{X} = \Max\left\{
\wgtBB{X;h} \midvert \text{$h$ is a generalised cohomology theory}
\right\}.
\end{align*}
\end{defn}
We often denote $\TCwgt{ \ ;h}{\pi}$ and 
$\wgtBB{ \ ;h}$ by $\TCwgt{ \ ;R}{\pi}$ and 
$\wgtBB{ \ ;R}$ respectively when $h^{\ast}( \ ) = H^{\ast}( \ ;R)$, where $R$ is a ring with unit.
%
%
%
%
%
We define versions of module weight for a fibrewise pointed space over $B$.
\begin{defn}\label{def:module-weight}
For a fibrewise pointed space $X$ over $B$, we define
\begin{enumerate}
\item
$\MwgtBB{X;h} = \Min\left\{m {\geq} 0%
\left\vert\,%
\text{\begin{minipage}{60mm}$(e^{X}_{m})^{\ast}$ is a split mono of (unstable) $h^{*}h$-modules\end{minipage}}
\right.\right\}$
for a generalisd cohomology theory $h$.
\item
$\MwgtBB{X} = \Max\left\{\MwgtBB{X;h} \midvert \text{$h$ is a generalised cohomology theory}\right\}$.
\end{enumerate}
\end{defn}

Then we immediately obtain the following result.
\begin{thm}
For any fibrewise pointed space $X$ over $B$, we have \par
\begin{center}
$\cuplenBB{X} \leq \wgtBB{X} \leq \MwgtBB{X} \leq \catBB{X} \leq \catlenBB{X} \leq \CatBB{X}$.
\end{center}
\end{thm}
By Lemma \ref{lem:milnor}, we have the following as a corollary of Theorem \ref{thm:main4}.
\begin{cor}
For any space $B$ having the homotopy type of a locally finite simplicial complex, we obtain 
\begin{center}
$\TCcuplen{B}{\pi} \leq \TCwgt{B}{\pi} \leq \MwgtBB{\double{B}} \leq \TC{B}{-}1$ 
$\leq \catlenBB{\double{B}} \leq \CatBB{\double{B}}$.
\end{center}
\end{cor}
\section{Higher Hopf invariants}
For any fibrewise pointed map $f : \SuspBB{V} \to X$ in $\TopBB$, we have its adjoint $\ad{f} : V \to \LoopBB{X}$ such that 
$$
e^{X}_{1}{\comp}\SuspBB{\ad{f}} = f : \SuspBB{V} \to X.
$$
If $\catBB{X} \leq m$, then there is a fibrewise pointed map $\sigma : X \to P^{m}_{B}\LoopBB{X}$ in $\TopBB$ such that 
$$
e^{X}_{1}{\comp}\sigma \simBB \id_{X} : X \to X.
$$
Hence both the fibrewise maps $e^{X}_{1}{\comp}(\sigma{\comp}f)$ and $e^{X}_{1}{\comp}\SuspBB{\ad{f}}$ are fibrewise pointed homotopic to $f$ in $\TopBB$.
Then we have 
$$
e^{X}_{1}{\comp}\{\SuspBB{\ad{f}}-(\sigma{\comp}f)\} \simBB *_{B},
$$
where $\simBB$ denotes the fibrewise pointed homotopy and $*_{B}$ denotes the fibrewise trivial map in $\TopBB$.
Thus there is a fibrewise pointed map $H^{\sigma}_{m}(f) : \SuspBB{V} \to E^{m+1}_{B}\LoopBB{X}$ such that 
$$
p^{\LoopBB{X}}_{m}{\comp}H^{\sigma}_{m}(f) \simBB \SuspBB{\ad{f}}-(\sigma{\comp}f).
$$
\begin{defn}
Let $X$ be of $\catBB{X} \leq m$, $m \geq 0$.
For $f : \SuspBB{V} \to X$, we define 
\begin{enumerate}
\item
$H^{B}_{m}(f) = \left\{H^{\sigma}_{m}(f) \midvert e^{X}_{1}{\comp}\sigma \simBB \id_{X} \right\} \subset [\SuspBB{V},X]$,
\item
$\mathcal{H}^{B}_{m}(f) = \left\{(\SuspBB{})^{\infty}_{*}H^{\sigma}_{m}(f) \midvert e^{X}_{1}{\comp}\sigma \simBB \id_{X} \right\} \subset \left\{\SuspBB{V},X\right\}^{B}_{B}$,
\end{enumerate}
where, for two fibrewise spaces $V$ and $W$, we denote by $\left\{V,W\right\}^{B}_{B}$ the homotopy set of fibrewise stable maps from $V$ to $W$.
\end{defn}

\appendix
\section{Fibrewise homotopy pull-backs and push-outs}

In this paper, we are using $A_{\infty}$ structures which is constructed using tools in $\TopB$ and $\TopBB$ --- especially, finite homotopy limits and colimits, in other words, fibrewise homotopy pull-backs and push-outs in $\TopB$ and $\TopBB$.
We show in this section that such constructions are possible even when a fibrewise space has some singular fibres.

First we consider the fibrewise homotopy pull-backs in $\TopBB$: let $X$, $Y$, $Z$ and $E$ be fibrewise spaces over $B$ and $p : E \to Z$ be a fibrewise fibration in $\TopB$.
For any fibrewise map $f : X \to Z$ in $\TopB$, there exists a pull-back $X \xleftarrow{f^{*}p} f^{*}E \xrightarrow{\hat{f}} E$ of $X \xrightarrow{f} Z \xleftarrow{p} E$ as 
$$
f^{*}E = \left\{(x,e) \in X{\times_{B}}E \midvert f(x)=p(e)\right\}
$$
a subspace of $X{\times_{B}}E$ together with fibrewise maps $f^{*}p : f^{*}E \to X$ and $\hat{f} : f^{*}E \to E$ given by restricting canonical projections:
$$
(f^{*}p)(x,e)=x,\quad \hat{f}(x,e)=e.
$$
\begin{thm}[Crabb-James \cite{CJ:fibrewise}]
Let $p : E \to Z$ be a fibrewise fibration.
For any fibrewise map $f : W \to Z$ in $\TopB$, $f^{*}p : f^{*}E \to W$ is also a fibrewise fibration.
\end{thm}
Let $\pi_{t} : \PathB{Z} \to Z$ be fibrewise fibrations given by $\pi_{t}(\ell)=\ell(t)$, $t=0,1$ (see also \cite{CJ:fibrewise}).
Then $\pi_{0}$ and $\pi_{1}$ induce a map $\pi : \PathB{Z} \to Z{\times_{B}}Z$ to the fibre product of two copies of $p_{Z} : Z \to B$.
\begin{prop}
$\pi : \PathB{Z} \to Z{\times_{B}}Z$ is a fibrewise fibration.
\end{prop}
\begin{proof}
For any fibrewise map $\phi : W \to \PathB{Z}$ and a fibrewise homotopy $H : W{\times}[0,1] = W{\times_{B}}(I_{B}) \to Z{\times_{B}}Z$ such that $H(w,0)=\pi{\comp}\phi(w)$ for $w \in W$, we define a fibrewise homotopy $\hat{H} : W{\times}[0,1] = W{\times_{B}}(I_{B}) \to \PathB{Z} (\subset \Path{Z})$ by 
$$
\hat{H}(w,s)(t) = \begin{cases}
\proj_{0}{\comp}H(w,s), & \text{if $t=0$,} 
\\
\proj_{0}{\comp}H(w,s{-}3t), & \text{if $0 < t < \frac{s}{3}$,} 
\\
\pi_{0}{\comp}\phi(w), & \text{if $t = \frac{s}{3}$,} 
\\
\phi(w)(\frac{3t-s}{3-2s}), & \text{if $\frac{s}{3} < t < \frac{3-s}{3}$,} 
\\
\pi_{1}{\comp}\phi(w), & \text{if $t = \frac{3-s}{3}$,} 
\\
\proj_{1}{\comp}H(w,3t{-}3{+}s), & \text{if $\frac{3-s}{3} < t < 1$}
\\
\proj_{1}{\comp}H(w,s), & \text{if $t=0$,} 
\end{cases}
$$
for $(w,s) \in W{\times_{B}}I_{B}$ and $t \in [0,1]$, where $\proj_{k} : Z{\times_{B}}Z \subset Z{\times}Z \to Z$ denotes the canonical projection given by $\proj_{k}(z_{0},z_{1})=z_{k}$,  $k=0,1$ for any $(z_{0},z_{1}) \in Z{\times_{B}}Z$.
Then for any $(w,s) \in W{\times_{B}}I_{B}$, we clearly have 
\begin{align*}&
\hat{H}(w,0)(t)=\phi(w)(t),\quad t \in [0,1],
\\&
(\hat{H}(w,s)(0),\hat{H}(w,s)(1)) = (\proj_{0}{\comp}H(w,s),\proj_{1}{\comp}H(w,s)) = H(w,s),
\end{align*}
and hence we have $\hat{H}(w,0)=\phi(w)$ for any $w \in W$ and also $\pi{\comp}\hat{H} = H$.
This implies that $\hat{H}$ is a fibrewise homotopy of $\phi$ covering $H$.
Thus $\pi$ is a fibrewise fibration.
\end{proof}
This yields the following corollary.
\begin{cor}
For any fibrewise maps $f : X \to Z$ and $g : Y \to Z$ in $\TopB$, the induced map $(f{\times_{B}}g)^{*}\pi : (f{\times_{B}}g)^{*}\PathB{Z} \to X{\times_{B}}Y$ is a fibrewise fibration in $\TopB$.
\end{cor}
We often call the fibrewise space $(f{\times_{B}}g)^{*}\PathB{Z}$ together with the projections $\proj_{X}{\comp}(f{\times_{B}}g)^{*}\pi : (f{\times_{B}}g)^{*}\PathB{Z} \to X$ and $\proj_{Y}{\comp}(f{\times_{B}}g)^{*}\pi : (f{\times_{B}}g)^{*}\PathB{Z} \to Y$ the homotopy pull-back in $\TopB$ of $X \xrightarrow{f} Z \xleftarrow{g} Y$.
We remark that the above construction can be performed within $\TopBB$ if $X$, $Y$, $Z$, $f$ and $g$ are all in $\TopBB$, so that we have a pointed version of a fibrewise homotopy pull-back:
\begin{cor}
For any fibrewise maps $f : X \to Z$ and $g : Y \to Z$ in $\TopBB$, the induced map $(f{\times_{B}}g)^{*}\pi : (f{\times_{B}}g)^{*}\PathB{Z} \to X{\times_{B}}Y$ is a fibrewise fibration in $\TopBB$.
\end{cor}

Second we consider the fibrewise homotopy push-outs in $\TopBB$: let $X$, $Y$, $Z$ and $W$ be fibrewise pointed spaces over $B$ and $i : Z \to W$ be a fibrewise cofibration in $\TopBB$.
For any fibrewise map $f : Z \to X$ over $B$, there exists a push-out $X \xrightarrow{f_{*}i} f_{*}W \xleftarrow{\check{f}} W$ of $X \xleftarrow{f} Z \xrightarrow{i} W$ as 
a quotient space of $X{\amalg_{B}}W$ by gluing $f(z)$ with $i(z)$ together with fibrewise maps $f_{*}i$ and $\check{f}$ induced from the canonical inclusions.
\begin{thm}[Crabb-James \cite{CJ:fibrewise}]
Let $i : Z \to W$ be a fibrewise cofibration in $\TopB$ (or $\TopBB$).
For any fibrewise map $f : Z \to X$ in $\TopB$ (or $\TopBB$, resp.), $f_{*}i : X \to f_{*}W$ is also a fibrewise cofibration in $\TopB$ (or $\TopBB$, resp.).
\end{thm}
Let us recall that $\TrackBB{Z}$ is obtained from $\TrackB{Z}=Z{\times_{B}}(B{\times}[0,1])=Z{\times}[0,1]$ by identifying the subspace $s_{Z}(B){\times}[0,1] \subset Z{\times}[0,1]$ with $s_{Z}(B)$ by the canonical projection to the first factor : $s_{Z}(B){\times}[0,1] \to s_{Z}(B)$.
Let $\iota_{t} : Z \to \TrackBB{Z}$ be fibrewise cofibration in $\TopBB$ given by $\iota_{t}(z)=q(z,t)$, $0 \leq t \leq 1$, where $q : Z{\times}[0,1] \to \TrackBB{Z}$ denotes the identification map.
Then $\iota_{0}$ and $\iota_{1}$ induce a map $\iota : Z{\vee_{B}}Z \to \TrackBB{Z}$ from $Z{\vee_{B}}Z$ the push-out of two copies of $s_{Z} : B \to Z$.
\begin{prop}
$\iota : Z{\vee_{B}}Z \to \TrackBB{Z}$ is a fibrewise cofibration.
\end{prop}
\begin{proof}
For any fibrewise map $\phi : \TrackBB{Z} \to W$ and a fibrewise homotopy $H : (Z{\vee_{B}}Z){\times}[0,1] = (Z{\vee_{B}}Z){\times_{B}}I_{B} \to W$ such that $H(z,0)=\phi{\comp}\iota(z)$ for $z \in Z{\vee_{B}}Z$, we define a fibrewise homotopy $\check{H} : \TrackBB{Z}{\times}[0,1] = \TrackBB{Z}{\times_{B}}(I_{B}) \to W$ by 
$$
\check{H}(q(z,t),s) = \begin{cases}
H(\incl_{0}(z),s{-}3t), & \text{if $0 \leq t < \frac{s}{3}$,}
\\
\phi{\comp}\iota_{0}(z), & \text{if $t = \frac{s}{3}$,} 
\\
\phi(q(z,\frac{3t-s}{3-2s})), & \text{if $\frac{s}{3} < t < \frac{3-s}{3}$,} 
\\
\phi{\comp}\iota_{1}(z), & \text{if $t = \frac{3-s}{3}$,} 
\\
H(\incl_{1}(z),3t{-}3{+}s), & \text{if $\frac{3-s}{3} < t \leq 1$}
\end{cases}
$$
for $(q(z,t),s) \in \TrackBB{Z}{\times_{B}}I_{B}$, where $\incl_{k} : Z \hookrightarrow Z{\vee_{B}}Z$, $k=0,1$ denote the canonical inclusion given by $\incl_{0}(z)=(z,*_{b})$ and $\incl_{1}(z)=(*_{b},z)$, $b=p_{Z}(z)$ for any $z \in Z$.
Then for any $(q(z,t),s) \in \TrackBB{Z}{\times_{B}}I_{B}$, we clearly have 
\begin{align*}&
\check{H}(q(z,t))(0)=\phi(q(z,t)),
\\&
\check{H}(q(z,0))(s) = H(\incl_{0}(z),s),\quad
\check{H}(q(z,1))(s) = H(\incl_{1}(z),s),
\end{align*}
and hence we have $\check{H}(q(z,t))(0)=\phi(q(z,t))$ for any $q(z,t) \in \TrackBB{Z}$ and also $\check{H}{\comp}(\iota{\times_{B}}1_{I_{B}}) = H$.
This implies that $\check{H}$ is a fibrewise homotopy of $\phi$ extending $H$.
Thus $\iota$ is a fibrewise cofibration.
\end{proof}
This yields the following corollary.
\begin{cor}
For any fibrewise maps $f : Z \to X$ and $g : Z \to Y$ in $\TopBB$, the induced map $(f{\vee_{B}}g)_{*}\iota : X{\vee_{B}}Y \to (f{\vee_{B}}g)^{*}\TrackBB{Z}$ is a fibrewise cofibration in $\TopBB$.
\end{cor}
We often call the fibrewise space $(f{\vee_{B}}g)^{*}\TrackBB{Z}$ together with the inclusions $(f{\vee_{B}}g)_{*}\iota{\comp}\incl_{X} : X \to (f{\vee_{B}}g)_{*}\TrackBB{Z}$ and $(f{\vee_{B}}g)_{*}\iota{\comp}\incl_{Y} : Y \to (f{\vee_{B}}g)_{*}\TrackBB{Z}$ as homotopy push-out in $\TopBB$ of $X \xleftarrow{f} Z \xrightarrow{g} Y$.

Quite similarly for a fibrewise space $Z$ in $\TopB$, we obtain a fibrewise cofibration $\hat\iota : Z \amalg Z = Z{\times}\{0\} \cup Z{\times}\{1\} \hookrightarrow Z{\times}[0,1] = \TrackB{Z}$.
Thus we have the following.
\begin{cor}
For any fibrewise maps $f : Z \to X$ and $g : Z \to Y$ in $\TopB$, the induced map $(f{\amalg}g)_{*}\hat\iota : X{\amalg}Y \to (f{\amalg}g)^{*}\TrackB{Z}$ is a fibrewise cofibration in $\TopB$.
\end{cor}
Thus we also have an unpointed version of a fibrewise homotopy push-out.

%
%


%
%
%


\begin{thebibliography}{99}
%
%
%
%
%
%
%
%
%
\bibitem{CJ:fibrewise}
M.~C.~Crabb. and I.~M.~James, ``Fibrewise Homotopy Theory'', Springer Monographs in Mathematics, Springer-Verlag London, Ltd., London, 1998. 
%
%
%
\bibitem{Farber:motion-planning}
M.~Farber, {\em Topological complexity of motion planning}, Discrete Comput. Geom. \textbf{29} (2003), 211--221.
%

\bibitem{Farber:robot-motion-planning}
M.~Farber, {\em Topology of robot motion planning}, ``Morse theoretic methods in nonlinear analysis and in symplectic topology'', 185--230, NATO Sci. Ser. II Math. Phys. Chem., 217, Springer, Dordrecht, 2006.
%
\bibitem{FG:smp}
M.~Farber and M.~Grant, {\em Symmetric Motion Planning}, Topology and robotics, 85--104, Contemp. Math., 438, Amer. Math. Soc., Providence, RI, 2007.
%
\bibitem{Fox:ls-cat}
R.~H.~Fox, {\em On the Lusternik-Schnirelmann category}, Ann. of Math. (2) \textbf{42}, (1941), 333--370.
%
\bibitem{Ganea:category}
T.~Ganea, {\em Lusternik-Schnirelmann category and strong category}, Illinois. J. Math. \textbf{11} (1967), 417--427.
%
%
%
%
%
%
%
%
%
%
%
%
%
%
%
%
%
%
%
%
%
\bibitem{Iwase:counter-ls}
N.~Iwase, {\em Ganea's conjecture on Lusternik-Schnirelmann category}, Bull. Lon. Math. Soc., \textbf{30} (1998), 623--634.
%
%
\bibitem{Iwase:counter-ls-m}
N.~Iwase, {\em $A_\infty$-method in Lusternik-Schnirelmann category}, Topology \textbf{41} (2002), 695--723.
%
\bibitem{Iwase:ls-cat-bundle}
N.~Iwase, {\em Lusternik-Schnirelmann category of a sphere-bundle over a sphere}, Topology \textbf{42} (2003), 701--713.
%
\bibitem{Iwase:cat-length}
N.~Iwase, {\em Categorical length, relative L-S category and higher Hopf invariants}, preprint.
%
\bibitem{IM:higher-associativity}
N.~Iwase and M.~Mimura, {\em Higher homotopy associativity}, Algebraic Topology, (Arcata CA 1986), Lect. Notes in Math. { \bf 1370}, Springer Verlag, Berlin (1989) 193--220.
%
%
%
\bibitem{IS:func_qfib}
N.~Iwase and M.~Sakai, {\em Functors on the category of quasi-fibrations}, Topology Appl. \textbf{155} (2008), 1403--1409.
%
%
%
%
\bibitem{JM:fibrewise-category}
I.~M.~James and J.~R.~Morris, {\em Fibrewise category}, Proc. Roy. Soc. Edinburgh. \textbf{119A} (1991), 177--190.
%
\bibitem{James:intro-fibrewise}
I.~M.~James, 
      {\em Introduction to fibrewise homotopy theory}, ``Handbook of algebraic topology'', 169--194, North Holland, Amsterdam, 1995.
%
\bibitem{James:ls-category-survey}
I.~M.~James, {\em Lusternik-Schnirelmann Category}, ``Handbook of algebraic topology'', 1293--1310, North Holland, Amsterdam, 1995.
%
%
%
%
%
%
%
%
\bibitem{Milnor:homotopy-type}
J.~Milnor, {\em On Spaces Having the Homotopy Type of a $CW$-Complex}, Trans. Amer. Math. Soc. \textbf{90} (1959), 272--280.
%
%
%
%
%
%
%
\bibitem{Rudyak:ls-cat_mfds2}
Y.~B.~Rudyak, {\em On category weight and its applications}, Topology \textbf{38} (1999), 37--55.
%
\bibitem{Rudyak:ls-cat_mfds3}Y.~B.~Rudyak, On analytical applications of stable homotopy (the Arnold conjecture, critical points), Math. Z. \textbf{230}(1999) 659--672.
%
%
\bibitem{Sakai:thesis}
M.~Sakai, {\em The functor on the category of quasi-fibrations}, DSc Thesis (Kyushu University 1999), 1999.
%
\bibitem{Sakai:push_pull}
M.~Sakai, {\em A proof of the homotopy push-out and pull-back lemma}, Proc. Amer. Math. Soc. \textbf{129} (2001), 2461--2466.
%
\bibitem{Schwarz:genus}
A.~S.~Schwarz {\em The genus of a fiber space}, Amer. Math. Soc. Transl.(2) \textbf{55} (1966), 49--140.
%
%
%
%
%
%
\bibitem{Stasheff:higher-associativity}
J.~D.~Stasheff, {\em Homotopy associativity of H-spaces, I, II}, Trans. Amer. Math. Soc. \textbf{108} (1963), 275--292, 293--312.
%
\bibitem{Strom:essential-cat-wgt}
J.~Strom, {\em Essential category weight and phantom maps}, Cohomological methods in homotopy theory (Bellaterra, 1998), 409--415, Progr. Math., 196, Birkhauser, Basel, 2001.
%
%
%
%
%
%
%
\bibitem{Whitehead:elements}
G.~W.~Whitehead, ``Elements of Homotopy Theory'', Springer Verlag, Berlin, GTM series \textbf{61}, 1978.
%
%
\end{thebibliography}
\end{document}


%
%
\begin{abstract}
There is a problem with the proof of Theorem 1.13 of \cite{IS10:tc-ls} which states that for a fibrewise well-pointed space $X$ over $B$, we have $\catBB{X} = \catBb{X}$ and that for a locally finite simplicial complex $B$, we have $\TC{B}=\TCM{B}$.
While we still conjecture that Theorem 1.13 is true, this problem means that, at present, no proof is given to exist.
Alternatively, we show the difference between two invariants $\catBb{X}$ and $\catBB{X}$ is at most $1$ and the conjecture is true for some cases.
We give further corrections mainly in the proof of Theorem 1.12.
\end{abstract}
%
%
\keywords{Toplogical complexity, Lusternik-Schnirelmann category.}
%
%
\subjclass[2000]{Primary 55M30, Secondary 55Q25}
%
%
\date{\today}
%
%
\maketitle
It was pointed out to the authors by Jose Calcines that there is a problem with the proof of Theorem 1.13 of \cite{IS10:tc-ls} which states that for a fibrewise well-pointed space $X$ over $B$, we have $\catBB{X} = \catBb{X}$ and that for a locally finite simplicial complex $B$, we have $\TC{B}=\TCM{B}$, where $\catBb{X}$ and $\TCM{B}$ are new versions of a fibrewise L-S category and a topological complexity, respectively, which are introduced in \cite{IS10:tc-ls}.

While we still conjecture that Theorem 1.13 of \cite{IS10:tc-ls} is true, this problem means that, at present, no proof is given to exist.
It then results that ``$\TC{B}$'' in Corollary 8.7 of \cite{IS10:tc-ls} must be replaced with ``$\TCM{B}$'' and the resulting inequality should be presented in the following form:
$$
\TCcuplen{B}{\pi} \leq \TCwgt{B}{\pi} \leq \MwgtBB{\double{B}} \leq \TCM{B}{-}1 
\leq \catlenBB{\double{B}} \leq \CatBB{\double{B}}.
$$

The problem in the argument occurs on page 14 where a homotopy 
$$
\hat{\Phi}_{i} : \hat{U}_{i} \times [0,1] \to \hat{X}
$$
is given, while the definition of $\hat{\Phi}_{i}$ apparently is not well-defined.
Alternatively, we show here the difference between two invariants $\catBb{X}$ and $\catBB{X}$ is at most $1$ and the conjecture is true for some cases.
\begin{thm}\label{thm:tc-ls-cor-1}
For a fibrewise \textit{well-pointed} space $X$ over $B$, we have 
$\catBb{X} \leq \catBB{X} \leq \catBb{X} + 1$
which implies that, for a locally finite simplicial complex $B$, we have 
$\TC{B} \leq \TCM{B} \leq \TC{B} + 1$.
\end{thm}
\begin{proof}
The inequality of $\TC{B}$ and $\TCM{B}$ in Theorem \ref{thm:tc-ls-cor-1} for a locally finite simplicial complex $B$ is, by Theorem 1.7 in \cite{IS10:tc-ls}, a special case of the inequality of $\catBb{X}$ and $\catBB{X}$ in Theorem \ref{thm:tc-ls-cor-1} for a fibrewise well-pointed space $X$.
So it is sufficient to show the inequality for $X$:
because the inequality $\catBb{X} \leq \catBB{X}$ is clear by definition, all we need to show is the inequality $\catBB{X} \leq \catBb{X}+1$.
Let $X$ be a fibrewise well-pointed space over $B$ with a projection $\proj_{X} : X \to B$ and a section $\sect_{X} : B \to X$.
Let $(u,h)$ be a fibrewise (strong) Str{\o}m structure (see Crabb and James \cite{CJ:fibrewise}) on $(X,\sect_{X}(B))$, i.e., $u : X \to [0,1]$ is a map and $h : X \times [0,1] \to X$ is a fibrewise pointed homotopy such that $u^{-1}(0) = \sect_{X}(B)$, $h(x,0)=x$ for any $x \in X$ and $h(x,1) = \sect_{X}\comp\proj_{X}(x)$ for any $x \in X$ with $u(x)<1$.
Assume $\catBb{X}=m$ and the family $\{U_{i}\,;\,0 {\leq} i {\leq} m\}$ of open sets of $X$ satisfies $X = \underset{i=0}{\overset{m}{\bigcup}} U_{i}$ and each open set $U_{i}$ is fibrewise contractible (into $\sect_{X}(B)$) by a fibrewise homotopy $H_{i} : U_{i} \times [0,1] \to X$.
Let $V_{i}=U'_{i} \cup V$ for $0 \leq i \leq m$ and $V_{m+1} = u^{-1}([0,\frac{2}{3}))$ where $U'_{i} = U_{i} \smallsetminus u^{-1}([0,\frac{1}{2}])$ and $V=u^{-1}([0,\frac{1}{3}))$.
Then the restriction $H_{i}\vert_{U'_{i}} : U'_{i} \times [0,1] \to X$ gives a fibrewise contraction of $U'_{i}$ and the restriction of the fibrewise (strong) Str{\o}m structure $h\vert_{V} : V \times [0,1] \to X$ gives a fibrewise pointed contraction of $V$.
Since $U'_{i}$ and $V$ are obviously disjoint, we obtain that $V_{i}=U'_{i} \cup V \supset \Delta(B)$ is a fibrewise contractible open set by a fibrewise pointed homotopy.
Similarly the restriction of the fibrewise (strong) Str{\o}m structure $h\vert_{V_{m+1}} : V_{m+1} \times [0,1] \to X$ gives a fibrewise pointed contraction of $V_{m+1} \supset \Delta(B)$.
Since $V_{i} \cup V_{m+1} = U'_{i} \cup V_{m+1} = U_{i} \cup V_{m+1} \supset U_{i}$, we obtain $\underset{i=0}{\overset{m+1}{\bigcup}} V_{i} = \underset{i=0}{\overset{m}{\bigcup}} \left(V_{i} \cup V_{m+1}\right) \supset \underset{i=0}{\overset{m}{\bigcup}} U_{i} = X$.
This implies $\catBB{X} \leq m+1 = \catBb{X} + 1$ and it completes the proof of Theorem \ref{thm:tc-ls-cor-1}.
\end{proof}
\begin{thm}\label{thm:tc-ls-cor-2}
Let $X$ be a fibrewise \textit{well-pointed} space over $B$ with $\catBb{X}=m$ and $\{U_{i}\,;\,0 {\leq} i {\leq} m\}$ be an open cover of $X$, in which $U_{i}$ is fibrewise contractible (into $\sect_{X}(B)$) by a fibrewise homotopy $H_{i} : U_{i} {\times} [0,1] \to X$.
Then we have $\catBB{X} = m = \catBb{X}$ if one of the following conditions is satisfied.
\begin{enumerate}
\item
There exists $i, \ 0 \leq i \leq m$ such that $U_{i}$ does not intersect with $\sect_{X}(B)$.
\item
There exists $i, \ 0 \leq i \leq m$ such that $U_{i}$ includes $\sect_{X}\comp\proj_{X}(U_{i}) \subset X$.
\end{enumerate}
\end{thm}
Theorem \ref{thm:tc-ls-cor-2} immediately implies the following corollary.
\begin{cor}\label{thm:tc-ls-cor-3}
Let $B$ be a locally finite simplicial complex with $\TC{B}=m$ and $\{U_{i}\,;\,1 {\leq} i {\leq} m\}$ be an open cover of $X$, in which $U_{i}$ is compressible into the image $\Delta(B)$ of diagonal map $\Delta : B \to B{\times}B$.
Then we have $\TCM{B} = m = \TC{B}$ if one of the following conditions is satisfied.
\begin{enumerate}
\item
There exists $i, \ 1 \leq i \leq m$ such that $U_{i}$ does not intersect with $\Delta(B)$.
\item 
There exists $i, \ 1 \leq i \leq m$ such that $U_{i}$ includes $\Delta\comp\pr_{2}(U_{i}) \subset B{\times}B$.
\end{enumerate}
\end{cor}
\begin{proof*}{Proof of Theorem \ref{thm:tc-ls-cor-2}:}
For simplicity, we assume that $i=0$ in each cases.
Let $X$ be a fibrewise well-pointed space over $B$ with a projection $\proj_{X} : X \to B$ and a section $\sect_{X} : B \to X$.
Let $(u,h)$ be a fibrewise (strong) Str{\o}m structure on $(X,\sect_{X}(B))$, i.e., $u : X \to [0,1]$ is a map and $h : X \times [0,1] \to X$ is a fibrewise pointed homotopy such that $u^{-1}(0) = \sect_{X}(B)$, $h(x,0)=x$ for any $x \in X$ and $h(x,1) = \sect_{X}\comp\proj_{X}(x)$ for any $x \in X$ with $u(x)<1$.
Then the fibrewise map $r : X \to X$ given by $r(x)=h(x,1)$ satisfies the following.
\begin{enumerate}
\item[i)]
$X = \underset{i=0}{\overset{m}{\bigcup}} r^{-1}(U_{i})$, since $X = \underset{i=0}{\overset{m}{\bigcup}} U_{i}$.
\item[ii)]
$r$ is fibrewise homotopic to the identity by $h$.
\item[iii)]
$r^{-1}(\sect_{X}(B)) \supset U = u^{-1}([0,1))$, where $U$ is fibrewise contractible by $h\vert_{U}$.
\item[iv)]
Each $r^{-1}(U_{i})$ is fibrewise contractible, since $r$ is fibrewise homotopic to the identity by ii) and $U_{i}$ is fibrewise contractible.
\end{enumerate}
Firstly, we consider the case (1): let $V_{0} = r^{-1}(U_{0}) \cup u^{-1}([0,\frac{2}{3}))$ and $V_{i} = (r^{-1}(U_{i}) \smallsetminus u^{-1}([0,\frac{1}{2}])) \cup u^{-1}([0,\frac{1}{3}))$, $1 \leq i \leq m$.
Thus $\underset{i=0}{\overset{m}{\bigcup}} V_{i} = r^{-1}(U_{0}) \cup \underset{i=1}{\overset{m}{\bigcup}} (V_{i} \cup u^{-1}([0,\frac{2}{3})))$ $\supset$ $r^{-1}(U_{0}) \cup \underset{i=1}{\overset{m}{\bigcup}} r^{-1}(U_{i})$ $=$ $\underset{i=0}{\overset{m}{\bigcup}} r^{-1}(U_{i})$ $=$ $X$ by i).
Since $r^{-1}(U_{i})$ is fibrewise contractible by iv), so is the open set $r^{-1}(U_{i}) \smallsetminus u^{-1}([0,\frac{1}{2}])$ for every $i \geq 0$, where $r^{-1}(U_{0}) \smallsetminus u^{-1}([0,\frac{1}{2}])=r^{-1}(U_{0})$ since $U_{0}$ does not intersect with $\sect_{X}(B)$.
On the other hand, $u^{-1}([0,\frac{t}{3}))$, $t=1,2$ are also fibrewise contractible by fibrewise pointed homotopies by iii).
Hence each $V_{i}$, $0 \leq i \leq m$ is fibrewise contractible by a fibrewise pointed homotopy, and hence $\catBB{X} \leq m = \catBb{X}$.
Thus we have $\catBb{X}=\catBB{X}$.

Secondly, we consider the case (2): 
let $V_{0} = r^{-1}(U_{0}) \cup u^{-1}([0,\frac{2}{3}))$ and $V_{i} = (r^{-1}(U_{i}) \smallsetminus u^{-1}([0,\frac{1}{2}])) \cup u^{-1}([0,\frac{1}{3}))$, $1 \leq i \leq m$.
Thus $\underset{i=0}{\overset{m}{\bigcup}} V_{i} = r^{-1}(U_{0}) \cup \underset{i=1}{\overset{m}{\bigcup}} (V_{i} \cup u^{-1}([0,\frac{2}{3})))$ $\supset$ $\underset{i=0}{\overset{m}{\bigcup}} r^{-1}(U_{i})$ $=$ $X$ by i).
Since $r^{-1}(U_{i})$ is fibrewise contractible by iv), so is the open set $r^{-1}(U_{i}) \smallsetminus u^{-1}([0,\frac{1}{2}])$ which does not intersect with $u^{-1}([0,\frac{1}{3}))$, for every $i>0$.
On the other hand, each open set $u^{-1}([0,\frac{t}{3}))$, $t=1,2$ is fibrewise contractible by a fibrewise pointed homotopy by iii).
Hence each open set $V_{i}$, $1 \leq i \leq m$ is fibrewise contractible by fibrewise pointed homotopy.
When $i=0$, we need to construct a fibrewise pointed homotopy $H : V_{0} \times [0,1] \to X$ using the fibrewise homotopy $H_{0} : U_{0} \times [0,1] \to X$ and the fibrewise (strong) Str{\o}m structure $(u,h)$ as follows:
$$
H(x,t) = \begin{cases}
\left\{\begin{array}{ll}
x,& t = 0
\\
h(x,3t),& 0 \leq t \leq \frac{1}{3}
\\
r(x),& t = \frac{1}{3}
\\
H_{0}(r(x),3t{-}1),& \frac{1}{3} \leq t \leq \frac{2}{3}
\\
\sect_{X}\comp\proj_{X}(r(x)) = \sect_{X}\comp\proj_{X}(x) = \sect_{X}(b),& t = \frac{2}{3}
\\
H_{0}(\sect_{X}(b),3{-}3t),& \frac{2}{3} \leq t \leq 1
\\
\sect_{X}(b),& t = 1
\end{array}\right\},& x \in V_{0} \smallsetminus U,
\\
\\[-4mm]
\left\{\begin{array}{ll}
x,& t = 0
\\
h(x,3t),& 0 \leq t \leq \frac{1}{3}
\\
r(x)=\sect_{X}(b),& t = \frac{1}{3}
\\
H_{0}(\sect_{X}(b),3t{-}1),& \frac{1}{3} \leq t \leq u(x){-}\frac{1}{3}
\\
H_{0}(\sect_{X}(b),3u(x){-}2),& u(x){-}\frac{1}{3} \leq t \leq \frac{5}{3}{-}u(x)
\\
H_{0}(\sect_{X}(b),3{-}3t),& \frac{5}{3}{-}u(x) \leq t \leq 1
\\
\sect_{X}(b),& t = 1
\end{array}\right\}, & \frac{2}{3} \leq u(x)<1,
\\
\\[-4mm]
\left\{\begin{array}{ll}
x,& t = 0
\\
h(x,3t),& 0 \leq t \leq \frac{1}{3}
\\
r(x)=\sect_{X}(b),& \frac{1}{3} \leq t \leq 1
\end{array}\right\}, & 0 \leq u(x) < \frac{2}{3},
\\
\\[-4mm]
\sect_{X}(b),& 
x \in \sect_{X}(B),
\end{cases}
$$
where $b = \proj_{X}(x) = \proj_{X}(r(x))$, and hence for $x \in V_{0} \smallsetminus u^{-1}([0,\frac{2}{3})) \subset r^{-1}(U_{0})$, we have $\sect_{X}(b) = \sect_{X}\comp\proj_{X}(r(x)) \in U_{0}$ since $r(x) \in U_{0}$.
Thus we have $\catBb{X}=\catBB{X}$, and it completes the proof of Theorem \ref{thm:tc-ls-cor-2}.
\end{proof*}
The following are corrections in \cite{IS10:tc-ls}.
\begin{itemize}
\item
The part of Proof of Theorem 1.12 from page 13 line -3 to page 14 line 12 is not clearly given and must be rewritten completely:\par
\begin{proof}
For each vertex $\beta$ of $B$, let $V_{\beta}$ be the star neighbourhood in $B$ and $V = \bigcup_{\beta} V_{\beta}{\times}V_{\beta} \subset B{\times}B$.
Then the closure $\bar{V} = \bigcup_{\beta} \bar{V}_{\beta}{\times}\bar{V}_{\beta}$ is a subcomplex of $B{\times}B$.
For the barycentric coordinates $\{\xi_{\beta}\}$ and $\{\eta_{\beta}\}$ of $x$ and $y$, resp, we see that $(x,y) \in V$ if and only if $\sum_{\beta}\Min(\xi_{\beta},\eta_{\beta})>0$ and that $\sum_{\beta}\Min(\xi_{\beta},\eta_{\beta})=1$ if and only if the barycentric coordinates of $x$ and $y$ are the same, or equivalently, $(x,y) \in \Delta(B)$.
Hence we can define a continuous map $v : B{\times}B \to [0,3]$ by the following formula.
\begin{align*}
v(x,y)=\begin{cases}
3-3\sum_{\beta}\Min(\xi_{\beta},\eta_{\beta}), & \text{if $(x,y) \in \bar{V}$,}\\
3, & \text{if $(x,y) \not\in V$.}
\end{cases}
\end{align*}
Since $B$ is locally finite, $v$ is well-defined on $B{\times}B$, and we have $v^{-1}(0) = \Delta(B)$ and $v^{-1}([0,3))=V$.
Let $U=v^{-1}([0,1))$ an open neighbourhood of $\Delta(B)$.
In \cite{Milnor:homotopy-type}, Milnor defined a map $\mu : V \to B$ giving an `average' of $(x,y) \in V$ as follows.\par
$$
\mu(x,y) = (\zeta_{\beta}),\quad \zeta_{\beta} = \Min(\xi_{\beta},\eta_{\beta})/\sum_{\gamma}\Min(\xi_{\gamma},\eta_{\gamma}),
$$
where $\{\xi_{\beta}\}$ and $\{\eta_{\beta}\}$ are barycentric coodinates of $x$ and $y$ respectively, and $\gamma$ runs over all vertices in $B$.
Since $B$ is locally finite, $\mu$ is well-defined on $V$ and satisfies $\mu(x,x)=x$ for any $x \in B$.
Using the map $\mu$, Milnor introduced a map $\lambda : V{\times}[0,1] \to B$ as follows.\par$$
\lambda(x,y,t) = 
\begin{cases}
(1{-}2t)x + 2t\mu(x,y),&t \leq \frac{1}{2},
\\
(2{-}2t)\mu(x,y) + (2t{-}1)y,&t \geq \frac{1}{2}.
\end{cases}
$$
Hence we have $\lambda(x,x,t)=x$ for any $x \in B$ and $t \in [0,1]$.
Using Milnor's map $\lambda$, we obtain a pair of maps $(u,h)$ as follows:
\begin{align*}&
u(x,y) = \Min\{1,v(x,y)\}\quad\text{and}
\\&
h(x,y,t) = \begin{cases}
(\lambda(x,y,\Min\{t,w(x,y)\}),y), 
& \text{if $v(x,y) < 3$,}\\
(x,y), & \text{if $v(x,y) > 2$,}
\end{cases}
\end{align*}
where $w : B{\times}B \to [0,1]$ is given by
$$
w(x,y)=
\begin{cases}
1,& v(x,y) \leq 1,
\\
2-v(x,y),&1\leq v(x,y)\leq2,
\\
0,&v(x,y) \geq 2.
\end{cases}
$$
If $2 < v(x,y) < 3$, then, by definition, we have $w(x,y)=0$ and 
$$
(\lambda(x,y,\Min\{t,w(x,y)\}),y) = (\lambda(x,y,0),y)=(x,y).
$$
Thus $h$ is also a well-defined continuous map.
Then we have $u^{-1}(0)=\Delta(B)$, $u^{-1}([0,1))=U$ and $h(x,y,0)=(x,y)$ for any $(x,y) \in B{\times}B$.
If $(x,y) \in U$, we have $w(x,y)=1$, $h(x,y,t) = (\lambda(x,y,t),y)$ and $h(x,y,1) = (y,y) \in \Delta(B)$.
Moreover, we have $h(x,x,t)=(x,x)$ for any $x \in B$ and $t \in [0,1]$ and $\pr_{2}\comp h(x,y,t) = y$ for any $(x,y,t) \in B{\times}B{\times}[0,1]$.
This implies that $h$ is a fibrewise pointed homotopy.
Thus the data $(u,h)$ gives the fibrewise (strong) Str{\o}m structure on $(B{\times}B,\Delta(B))$.
\end{proof}
\item
In page 19, line 34, ``$t=0$'' must be replaced by ``$t=1$''.
\item
In page 20, line 17, ``that'' must be replaced by ``that $H(s_{Z}(b),t)=s_{W}(b)$ for any $b \in B$ and''.
\item
In page 20, line 28, the formula ``$\check{H}(q(s_{Z}(b),t),s)=s_{W}(b)$,'' must be added.
\end{itemize}